\documentclass[12pt]{article}
\usepackage{latexsym}
\usepackage{graphicx}
\usepackage{verbatim}
\usepackage{amssymb}
\usepackage{amsmath}
\title{A  H\"{o}lder-type inequality on a regular rooted tree} 
\author{K.J. Falconer\\
\small{{\it Mathematical Institute,  
University of St~Andrews, North Haugh, St~Andrews,}} \\
\small{{\it Fife, KY16~9SS, Scotland }}} 
\date{}

\newtheorem{theo}{Theorem}

\newtheorem{lem}[theo]{Lemma}
\newtheorem{cor}[theo]{Corollary}

\newcommand{\bi}{{\bf i}} 
\newcommand{\bj}{{\bf j}} 
\newcommand{\bv}{{\bf v}} 
\newcommand{\bw}{{\bf w}} 
 
\newcommand{\be}{\begin{equation}} 
\newcommand{\ee}{\end{equation}} 
\allowdisplaybreaks
\renewcommand{\a}{\alpha}

\begin{document}
\maketitle

\begin{abstract} 
We establish an inequality which involves a non-negative function defined on the vertices of a finite $m$-ary regular rooted tree. The inequality may be thought of as relating an interaction energy defined on the free vertices of the tree summed over automorphisms of the tree, to a product of sums of powers of the function over vertices at certain  levels of the tree. Conjugate powers arise naturally in the inequality, indeed, H\"{o}lder's inequality is a key tool in the proof which uses induction on subgroups of the automorphism group of the tree.
\end{abstract}

\medskip

 \section{Introduction}
\setcounter{equation}{0}
\setcounter{theo}{0}

The energy of an interacting particle system with a finite number of sites 
is typically given by a sum $\sum_{(\bi_1, \bi_2)} \mu(\bi_1)\mu(\bi_2)F(\bi_1,\bi_2)$ over pairs of particles $(\bi_1, \bi_2)$  where $F(\bi_1,\bi_2)$ represents the force or interaction between the particles located at $\bi_1$ and $\bi_2$ which have masses or weights $\mu(\bi_1)$ and $\mu(\bi_2)$ respectively.

A more complex system may involve interactions between three or more particles simultaneously rather than just two. Thus we are led to consider energies of the form 
$\sum_{(\bi_1, \ldots, \bi_n)} \mu(\bi_1)\cdots \mu(\bi_n)F(\bi_1,\ldots, \bi_n)$ resulting from interactions $F(\bi_1,\ldots, \bi_n)$ which depend on $n$ particles and their configuration. Typically each particle may be affected most by those other particles that are closest, so the interaction should take into account the nearest neighbor structure of the particle configuration. One convenient way of incorporating such a structure is by representing the sites as the free vertices (i.e. vertices of valence 1) of a regular $m$-ary tree, so that the distance between a pair of sites is an ultrametric determined by the level of their first common ancestor. The arrangement of common ancestors or `joins' of a collection of $n$ particles determines their nearest neighbor configuration.

Our main results and detailed notation are set out in Section 2, but to fix ideas we give here a brief overview and an example.  We work on an $m$-ary regular rooted tree $T$ of $k$ levels, where the level or  generation of a vertex is its edge distance from the root. In the usual parlance, each vertex has $m$ `children', except for  the free vertices at level $k$ which have no children; we let $T_0$ denote the set of free vertices. We write $\bi \wedge \bi' $ for  the {\it join} of $\bi, \bi' \in T_0$, that is the vertex $\bj \in T$ of maximum level that lies on both of the paths from the root $\emptyset$ to $\bi$ and from $\emptyset$ to $ \bi' $.   The  {\it join set} of  $n$ particles located at sites $\bi_1,\ldots,\bi_n \in T_0$,  is the set of vertices of $T$ given by $\bi_i\wedge \bi_j $ for all $1\leq i < j\leq n$ (this is made a little more precise in the next section where we allow for multiple join points). 

For $f$ a non-negative function defined on the vertices of $T$, we will consider interactions that can be expressed as
the product of the values  of $f$ over the join points, that is,
$$F(\bi_1,\ldots,\bi_n) = f(\bj_1)f(\bj_2) \cdots f(\bj_{n-1})$$
 where $\bj_1,\bj_2, \ldots, \bj_{n-1}$ are the join points of $\bi_1,\ldots,\bi_n$.
  Then the energy is of the form
\begin{equation}
\sum_{(\bi_1,\ldots,\bi_n)} \mu(\bi_1)\cdots\mu(\bi_n)F(\bi_1,\ldots,\bi_{n}).\label{ensum}
\end{equation}
Typically, $f(\bj)$ will be large if the join point $\bj$ is a high level vertex corresponding to a significant interaction component resulting from nearby particles. For many problems one needs to bound the energy of a system, perhaps in the limit as the number of generations becomes large. Such estimates are required, for example, in estimating high moments of certain measures, see for \cite{Fa5,FX}.

Whilst one may wish to estimate the sum (\ref{ensum}) over all arrangements of $n$ particles, the sum breaks down naturally into sub-sums over configurations of  particles which have isomorphic join structures. Thus we consider the set of configurations of $n$ particles on $T_0$ that may be obtained from each other under some automorphism of the rooted tree, in other words the equivalence classes of configurations defined by the automorphisms. These equivalence classes are the orbits of the automorphism group of T acting on the ordered $n$-tuples from $T_0$. Writing $[I]$ for such an equivalence class or orbit, we are led to consider the sums
$$\sum_{(\bi_1,\ldots,\bi_n)\in [I]} \mu(\bi_1)\cdots\mu(\bi_n)F(\bi_1,\ldots,\bi_n).$$
We will obtain upper bounds for these sums over the equivalence class $[I]$ in terms of $p$th  powers of   $f$ summed across certain levels of the tree $T$.
\begin{figure}[h]\label{Fig1}
\begin{center}
\includegraphics[scale=0.5,bb=50 250 500 600]{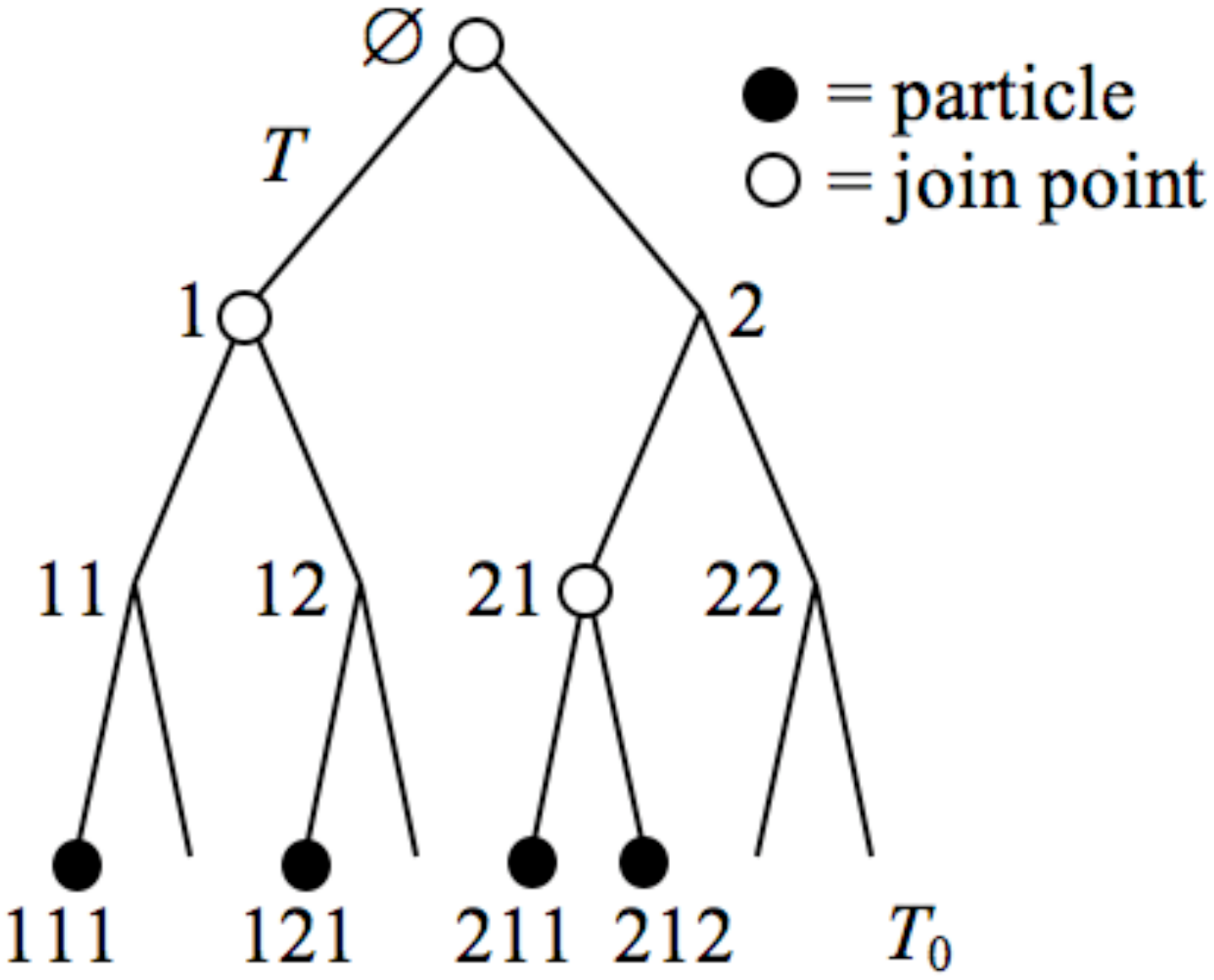}
\caption{Particles and join points in the example}
\end{center}
\end{figure}
 
To illustrate this, consider a specific case of the binary rooted tree of  $k=3$ levels with $n=4$ particles. Starting with the configuration  $(111,121,211,212)$ of points of $T_0$ (with the usual coding of vertices of the binary tree) the join points are at $\emptyset, 1$ and $21$, see Figure 1. Let $[I]$ be the class of 64 different equivalent ordered 4-tuples  $(g(111),g(121),g(211),g(212))$  obtainable under an automorphism $g$ of the rooted tree $T$. Assign each free vertex $\bi \in T_0$ a positive weight $\mu(\bi)$. For each vertex $\bj \in T$ write $\mu(\bj)$ for the total weight of the free vertices below $\bj$, so, for example, $\mu(21) = \mu(211)+\mu(212)$.
In this special case, our main inequality becomes,  for any non-negative function $f$ on the vertices of  $T$ and for all  $p_1,p_2,p_3>2 $  satisfying 
$\frac{1}{p_1} +\frac{1}{p_2} +\frac{1}{p_3} =1$,
\begin{align}
\sum_{(\bi_1,\bi_2,\bi_3,\bi_4) \in[I]}  \mu(\bi_1)\mu(\bi_2)&\mu(\bi_3)\mu(\bi_4)F(\bi_1,\bi_2,\bi_3,\bi_4)\label{eqn1}\\
&\leq \, 8\cdot 8^{1/p_1} \cdot4^{1/p_2}\cdot 2^{1/p_3}
\big(f(\emptyset)\mu(\emptyset)^{1+p_1}\big)^{1/p_1}\nonumber\\
&\quad\times \big(f(1)^{p_2}\mu(1)^{1+p_2} +f(2)^{p_2}\mu(2)^{1+p_2}\big)^{1/p_2}\nonumber\\
&
\quad\times \big(f(11)^{p_3}\mu(11)^{1+p_3} +f(12)^{p_3} \mu(12)^{1+p_3}\nonumber\\
&\qquad\qquad \qquad+f(21)^{p_3}\mu(21)^{1+p_3}  +f(22)^{p_3}\mu(22)^{1+p_3}\big)^{1/p_3}. \label{eqn2}
\end{align}
The sum (\ref{eqn1}) has 64 terms, each a product over  3 join points, one at each level 0, 1 and 2. The bound (\ref{eqn2}) is a product of weighted $p_i$th power sums of $f$ across each of the 3 levels of the tree at which the join points occur.
A check shows that equality holds if the weights $\mu(\bi)$ are equal for all $\bi \in T_0$ and $f$ is constant across each level, that is if  $f(1) = f(2)$ and $f(11) =f(12)= f(21) =f(22)$.

We will prove results of this type in a general context by employing induction with respect to automorphisms fixing the  join points and incorporating H\"{o}lder's inequality in a natural way. 

\section{Notation and statement of results}
\setcounter{equation}{0}
\setcounter{theo}{0}

To  state the main results we set out some further notation relating to the tree $T$.

For integers $m\geq 2$ and $k \geq1$, we index the vertices of $T$, the $m$-ary regular rooted tree of $k+1$ levels (including the level of the root), by the {\it symbolic space} of words formed from the symbols $\{1,2,\ldots,m\}$. 
Thus  the vertices of $T$ are given by $\{ (i_{1}, i_{2}, \ldots , i_{l}): 0 \leq l\leq k, \, 1 \leq i_{j} \leq 
m \}$, with the root of the tree as the empty word $\emptyset$. 
 We often abbreviate a
word by
${\bf i} = (i_{1}, i_{2}, \ldots , i_{l} ) $ and write $|{\bf i}|=l$ for 
its {\it level} or {\it generation}. The set of  free vertices, that is those $\bi$ with $|\bi|= k$, is denoted by $T_0$.
We write ${\bf i} \preceq {\bf j}$ to mean that ${\bf i}$ is a 
{\it curtailment} of ${\bf j}$, that is  ${\bf i}$ is an initial subword of ${\bf j}$.
If $\bi, \bi' \in T$  then the {\it join} 
${\bi}\wedge {\bi'}$ is the maximal word such that both 
${\bi}\wedge {\bi'} \preceq {\bi}$ and ${\bi}\wedge {\bi'} \preceq {\bi'}$.  With  each $\bj \in T$ we associate
the {\it cylinder}  
$C_\bj = \{\bi \in T_0 : \bj \preceq \bi \}$ comprising those points of $T_0$ below $\bj$.

We next consider automorphisms of the rooted tree $T$, regarded as a graph, which induce permutations of the vertices at each level of the tree.
 For each vertex $\bv \in T$  and $1 \leq n \leq m^k$ we write
\be 
S_\bv(n) = \big\{(\bi_1,\ldots,\bi_n) : \bi_j \in T_0, \, \bi_j  \succeq \bv \, (1 \leq j \leq n), \, \bi_j \neq \bi_h \, (j \neq h)\big\} \label{sdef}
\ee
for the set of all ordered $n$-tuples of distinct elements of $T_0$   that are descendents of  $\bv$. 
Let $\mbox{Aut}_\bv$ be the group of automorphisms of the rooted tree $T$ that fix $\bv$. Define an equivalence relation $\sim$ on $S_\bv(n)$
by 
\be 
(\bi_1,\ldots,\bi_n)  \sim (\bi_1',\ldots,\bi_n') \mbox{ if there exists } g \in   \mbox{Aut}_\bv  \mbox{ such that  }   g(\bi_r) =  \bi_r'  \mbox{ for  all }  1 \leq r \leq n; \label{equivdef}
\ee
  thus the  equivalence classes  are the orbits of $S_\bv(n)$ under $ \mbox{Aut}_\bv $. We write $[I]_\bv$ for the equivalence class containing $I = (\bi_1,\ldots,\bi_n)$.
For notational simplicity, we often omit the subscript when $\bv = \emptyset$, so that
$S(n) = S_\emptyset(n)$ and $[I]= [I]_\emptyset$.

We require some terminology relating to the joins of  subsets  of  $T_0$. Let $I =  (\bi_1,\ldots,\bi_n)\in S_\bv(n)$. The {\it join set} of $I$, denoted by $\bigwedge(I) = \bigwedge(\bi_1,\ldots,\bi_n)$,  is the set of vertices $\{\bj_1,\ldots,\bj_{n-1}\} \subset T$ consisting of the {\it join points}  $\bi_i\wedge \bi_j $ for all $1\leq i <j \leq n$, with $\bw \in \bigwedge(I)$ occurring with {\it multiplicity} $r$ if there are $(r+1)$ distinct indices $1 \leq i_1< \ldots<i_{r+1}\leq n $ such that  
$\bi_{i_s}\wedge\bi_{i_t} = \bw$ for all $s \neq t$.  Note that the join set of $n$ points always consists of $n-1$ points counting by multiplicity. Moreover, if $T$ is a binary tree, so $m=2$, all join points have multiplicity $1$.

\begin{figure}[h]\label{Fig2}
\begin{center}
\includegraphics[scale=0.45,bb=160 180 500 650]{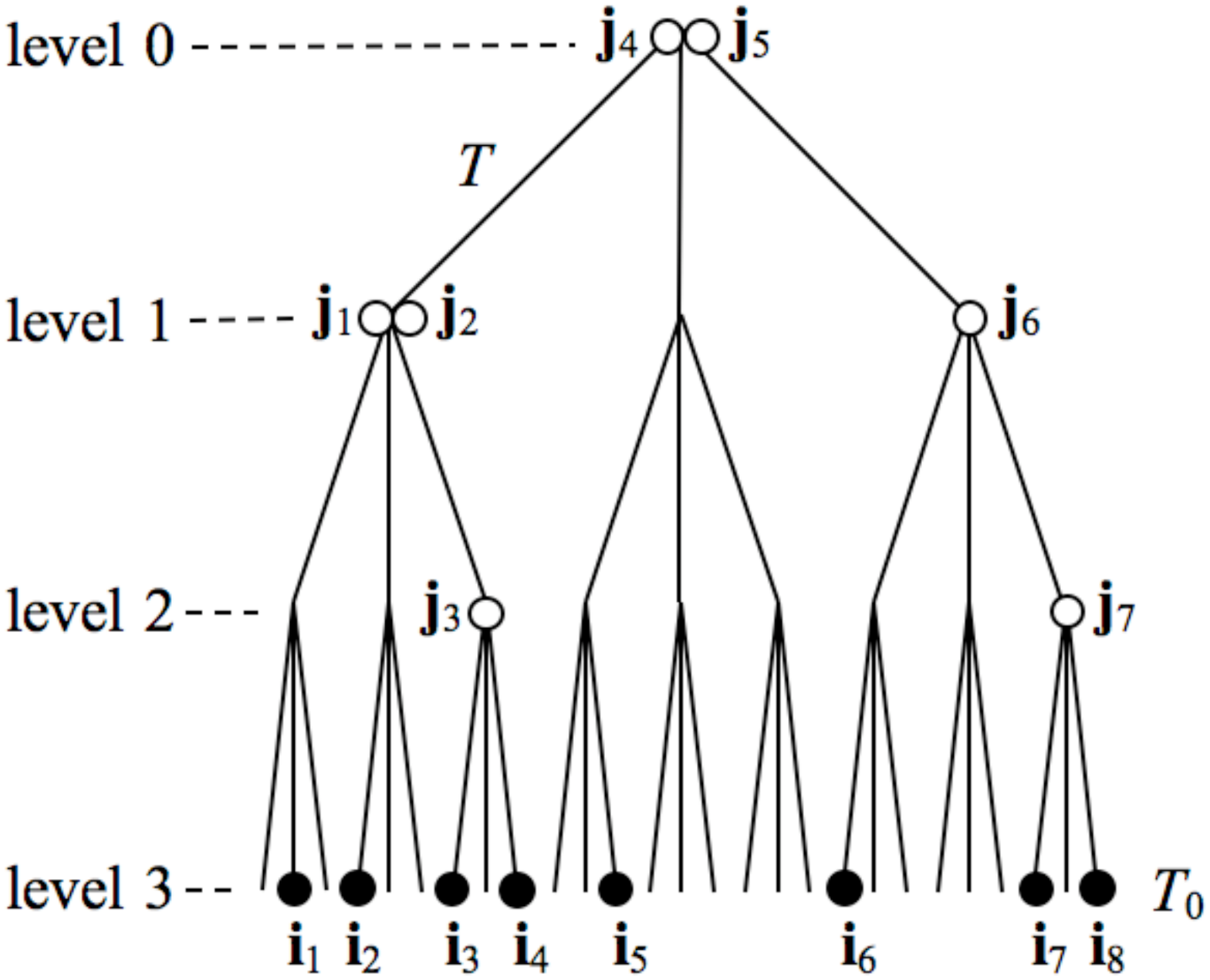}
\caption{A 3-ary tree with 8 particles and 7 join points, 2 of which have multiplicity 2}
\end{center}
\end{figure}

The    {\it set of join levels} $L(I)$ of $I \in S_\bv(n)$ is  $\{|\bj_1|,\ldots,|\bj_{n-1}|: \bj_i \in \bigwedge(I)\}$ with levels repeated according to multiplicity. Notice that if $I \sim I'$ then $L(I) = L(I') =: L([I]_\bv)$, i.e. the set of levels is constant across each equivalence class $[I]_\bv$ of $S_\bv(n)$.

For a level $L$ of $T$ we will, by slight abuse of notation, write $\bj \in L$ to mean that $|\bj| = L$, so we think of a vertex $\bj$ `belonging' to a level of the tree.

We now assign a weight $\mu(\bi)  \geq 0$ to each $\bi \in T_0$. Then  the weight of each cylinder is defined to be the sum of the weights of the points in the cylinder, that is   $\mu(C_\bj )= \sum_{\bi \in T_0, \bj \preceq \bi} \mu(\bi)$ for each $\bj \in T$.

Next we assign a positive value $f(\bj)$ to each vertex $\bj$ of $T$, that is  $f: T\to \mathbb{R}^+$.
Then for each $I \in S(n)$ we take the product of these values over the join points of $I$. Thus if
$n \geq 2$  we define $F: S(n) \to \mathbb{R}^+$ by
\be
F(\bi_1,\ldots,\bi_n) = f(\bj_1)f(\bj_2) \cdots f(\bj_{n-1}), \label{deff}
\ee
where $I =  (\bi_1,\ldots,\bi_n)\in S(n)$ has join set  $\bigwedge(\bi_1,\ldots,\bi_n) = \{\bj_1, \ldots, \bj_{n-1}\}$. For $I  = (\bi_1) \in S(1)$ we make the trivial assignment $F(\bi_1) =1$. Note that if $\bj$ is a join point of multiplicity $r$ then the factor $f(\bj)$ will occur in (\ref{deff}) $r$ times

We can now state the main theorem.

\begin{theo}\label{thmmary}
Let $T$ be the rooted $m$-ary tree  of  $k$ levels. Let   $2 \leq n \leq m^{k}$. Let  $I \in S(n)$ have join points at levels $L_1, \ldots,L_{n-1}$ and let
$p_1, \ldots,p_{n-1}>0$ satisfy $\sum_{i=1}^{n-1} \frac{1}{p_i} = 1$. Then
\be
\sum_{(\bi_1,\ldots,\bi_n) \in[I]} \mu(\bi_1)\cdots\mu(\bi_n)F(\bi_1,\ldots,\bi_n)
\leq K \prod_{i=1}^{n-1}
 \Big(\sum_{\bj\in L_i} f(\bj)^{p_i}\mu(C_{\bj})^{1+p_i}\Big)^{1/p_i}, \label{intesta}
\ee
 where $K$ does not depend on $f$ or $\mu$.
\end{theo}

The value of the constant $K$ will be discussed in Section \ref{valk}. We obtain the following estimate, though this is not in general optimal.

\begin{cor}\label{genk}
We may take $\displaystyle{K = \prod_{i=1}^j \frac{(m-1)! }{(m-r_i-1)! } } \leq (m-1)^{n-1}$ in \eqref{intesta}, where the product is over the distinct join points which have multiplicities $r_1,\ldots,r_j$. In particular, for a binary tree, where $m=2$, we can take $K=1$.
\end{cor}

For a binary tree  we can, under certain conditions, improve on this to obtain the best possible value of $K$. 

\begin{cor}\label{thmbinary}
Let $T$ be the rooted binary tree  of $k$ levels. Let   $2 \leq n \leq 2^{k}$. Let  $I \in S(n)$ have join points $\{\bj_1,\ldots,\bj_{n-1}\}$ at levels $L_1, \ldots,L_{n-1}$ and let
$p_1, \ldots,p_{n-1}>0$ satisfy $\sum_{i=1}^{n-1} \frac{1}{p_i} = 1$. Suppose also that the following condition is satisfied:

If the join points of $I$ are ordered so that $\bj_1$ is the `top' join point (i.e. $\bj_1 \preceq \bj_i$ for all $i$) and  
$\{\bj_2,\ldots,\bj_{m}\}$ and $\{\bj_{m+1},\ldots,\bj_{n}\}$ are the join points below the left and right edges abutting  $\bj_1$ respectively, then
\begin{equation}
\sum_{i=2}^{m} \frac{1}{p_i} \leq \frac{1}{2} \quad \mbox{ and } \sum_{i=m+1}^{n} \frac{1}{p_i} \leq \frac{1}{2}. \label{halves}
\end{equation}
Then
\be
\sum_{(\bi_1,\ldots,\bi_n) \in[I]} \mu(\bi_1)\cdots\mu(\bi_n)F(\bi_1,\ldots,\bi_n)
\leq \frac{1}{2^{n-1}} \prod_{i=1}^{n-1}
 \Big(\sum_{\bj\in L_i} f(\bj)^{p_i}\mu(C_{\bj})^{1+p_i}\Big)^{1/p_i}, \label{intestabin}
\ee
with equality if $\mu(\bi)$ is constant for all $\bi \in T_0$ and, for each level $L_i$, $f(\bj)$ is constant for all $\bj \in L_i$.
\end{cor}

\section{Proof of Theorem \ref{thmmary}}
\setcounter{equation}{0}
\setcounter{theo}{0}

We will establish  the key Lemma \ref{lemA} by induction on the join points, using repeated applications of 
H\"{o}lder's inequality and at each step `carrying forward' some of the `weight' to the next join point. 
Note that the proof is simpler when $T$ is a binary tree or when all the join points have multiplicity 1, in which case $d=2$ throughout the induction.
Theorem \ref{thmbinary}  follows immediately from the lemma on taking $\bv= \emptyset$ to be the root of the tree $T$.

In the proof of Lemma \ref{lemA} we will refer to the numbers $K(m;a_1, \ldots,a_m)$, defined as follows for $m \geq 1$ and $a_i \geq 0$. With  $S_m$ as the symmetric group of all permutations $\sigma$ of $(1,2,\ldots,m)$, $K(m;a_1, \ldots,a_m)$  is  the least number such that  
\be
  \sum_{\sigma \in S_m} x_{\sigma(1)}^{a_1}x_{\sigma(2)}^{a_2}\cdots x_{\sigma(m)}^{a_m}
  \leq K(m;a_1, \ldots,a_m)(x_1  + \cdots +x_m)^{a_1 + \cdots +a_m} \label{muirin1}
\ee
for all $x_i \geq0 $, where we make the consistent convention that $0^0 =1$. It is easy to see that $K(m;a_1, \ldots,a_m) \leq  m!$ ; we will discuss the value of $K$ in more detail in Section \ref{valk}.

\begin{lem}\label{lemA}
Let $T$ be the rooted $m$-ary tree  of  $k$ levels. Let  $\bv \in T$ and $1 \leq n \leq m^{k-|\bv|}$. Let  $I \in S_\bv(n)$ have join points at levels $L_1, \ldots,L_{n-1}$ and let
$p_1, \ldots,p_{n-1}>0$ satisfy $\sum_{i=1}^{n-1} \frac{1}{p_i} = 1$. Then
\be
\sum_{(\bi_1,\ldots,\bi_n) \in[I]_\bv} \mu(\bi_1)\cdots\mu(\bi_n)F(\bi_1,\ldots,\bi_n)
\leq K \prod_{i=1}^{n-1}
 \Big(\sum_{\bj\in L_i, \bj \succeq \bv} f(\bj)^{p_i}\mu(C_{\bj})^{1+p_i}\Big)^{1/p_i}, \label{intestapf}
\ee
where $K$ does not depend on $f$ or $\mu$.
\end{lem}

\noindent{\it Proof.} 
We proceed by strong induction. We take the inductive hypothesis \textbf{P}(\textit{n}) for $n \in \mathbb{N}$ to be as follows.
\medskip

\noindent\textbf{P}(\textit{n}) $(n \geq 1)$:
For  all $\bv \in T$, for all $p_1, \ldots,p_{n-1}>0$ and $1 \leq \alpha \leq \infty$ satisfying $\sum_{i=1}^{n-1} \frac{1}{p_i} = 1-\frac{1}{\alpha}$, if
$I \in S_\bv(n)$ has join levels $L_1, \ldots,L_{n-1}$,  there is a constant $K$ such that
\be
\sum_{(\bi_1,\ldots,\bi_n) \in[I]_\bv} \mu(\bi_1)\cdots\mu(\bi_n)F(\bi_1,\ldots,\bi_n)
\leq K \mu(C_{\bv})^{1/\a}\prod_{i=1}^{n-1}
 \Big(\sum_{\bj\in L_i, \bj \succeq \bv} f(\bj)^{p_i}\mu(C_{\bj})^{1+p_i}\Big)^{1/p_i}, \label{intestapnp}
\ee
for all  $\mu$ and $f$.

When $n=1$ then $\a = 1$ and we take the empty product to equal 1 in (\ref{intestapnp}).

The proof falls into several parts.
\medskip

\noindent (a)  {\it Verification of }{\bf P}(1).
If   $\bv \in  T$ and $I = (\bi_1) \in S_\bv(1)$ then  
$$
\sum_{(\bi_1) \in[I]_{\bv}} \mu(\bi_1) F(\bi_1) = \sum_{\bi_1 \succeq \bv} \mu(\bi_1)  = \mu(C_{\bv}),$$
which is $\textbf{P}({1})$ with $K=1$, noting that $\alpha=1$.
\medskip

\noindent (b) {\it The inductive step}.
Assume, for some integer $n \geq 2$, that  \textbf{P}(\textit{l})  holds for all  $1 \leq l\leq n-1$. We will show that \textbf{P}(\textit{n}) holds.
Let $\bv \in T$ and $I = (\bi_1,\ldots,\bi_n) \in S_\bv(n)$ where, without loss of generality, $(\bi_1,\ldots,\bi_n)$ are indexed in lexicographical order.  We consider two cases.
\smallskip

\noindent{\it Case} (b1). Suppose that $\bv \in T$ is a join point of $I$ of multiplicity $d-1 \ (2 \leq d \leq m)$. Let $\bv_1,\ldots, \bv_d$ be the vertices of $T$ that are children of $\bv$ and are above at least one of the $\bi_j$, that is for each $j$, $|\bv_j| = |\bv| +1$ and there is an $k$ such that 
$\bv \preceq \bv_j \preceq \bi_k$.
Then $I$ splits into $d$ non-empty parts,  $I_j= (\bi_1^j,\ldots,\bi_{n_j}^j)$  for $1 \leq j \leq d$, so that $\bi_i^j \succeq \bv_j$ for all $i,j$, with $n_j \geq 1$ and  $n_1 + \cdots + n_d = n$. Note that the join points of $I$ comprise the join points of each of the $I_j$ together with $\bv$.

Let $p_1, \ldots,p_{n-1}>0$  be the exponents associated with the join points of $I$ where $\sum_{i=1}^{n-1} \frac{1}{p_i} = 1- \frac{1}{\alpha}$ with $1 <\alpha \leq \infty$. Re-index these exponents so that  $p_1^j, \ldots, p_{n_j-1}^j$ are those associated with the join points of $I_j$ for each $j$ and the exponents  $p_1^0, \ldots, p_{d-1}^0$ are those associated with $\bv$. Let the corresponding join levels of $I_j$ be 
$\Lambda(I_j) = \{L_1^j,\ldots,L_{n_j -1}^j\}$ for $1 \leq j \leq d$.
Define $\a_1, \ldots, \a_d$  by
\begin{equation}
\sum_{i=1}^{n_j-1} \frac{1}{p_i^j} =  1-\frac{1}{\a_j}  < 1 \mbox{ (if $n_j\geq 2$) or }
 \a_j = 1\mbox{ (if $n_j=1$)}.\label{ajs}
 \end{equation}
Summing over $j$ gives
\begin{equation}
   \frac{1}{\a_1} + \cdots + \frac{1}{\a_d} = d-1+  \frac{1}{p_1^0} + \cdots + \frac{1}{p_{d-1}^0}+ \frac{1}{\a} = d-1 + \frac{1}{\beta}\label{betadef}
   \end{equation}
where, for convenience, we write $1/\beta = 1/p_1^0 + \cdots + 1/p_{d-1}^0 + 1/\a$.

Now let $(\bv_1,\ldots,\bv_d, \ldots, \bv_m)$ be the set of {\it all} vertices of $T$ immediately below $\bv$.   Let $G$ be a subgroup of $ \mbox{Aut}_\bv$  such that $\{(g(\bv_1),\ldots,g(\bv_m)): g \in G\}$ includes every permutation of   $(\bv_1,\ldots,\bv_m)$ exactly once, so that $G$ has order $m!$. The orbits of $I$ may be decomposed as
$$[I]_\bv = \{( [g(I_1)]_{g(\bv_1)},\ldots, [g(I_d)]_{g(\bv_d)})\}_{g \in G},$$
but if $m \geq d+2$ then $(m-d)!$ elements  $g$ of $G$  will yield each such decomposition of $[I]_\bv$.

Then, by definition of $F$,
\begin{align}
&(m-d)!\sum_{(\bi_1,\ldots,\bi_n) \in[I]_\bv}\mu(\bi_1)\cdots\mu(\bi_n)F(\bi_1,\ldots,\bi_n)\nonumber
\\
&=
\sum_{g\in G} f(\bv)^{d-1}\bigg(\sum_{(\bi_1^1,\ldots,\bi_{n_1}^1) \in[g(I_1)]_{g(\bv_1)}}
\!\!\!\!\mu(\bi_1^1)\cdots\mu(\bi_{n_1}^1)F(\bi_1^1,\ldots,\bi_{n_1}^1)\bigg)\nonumber
\\
&\hspace{4cm}\times \cdots \times\bigg(\sum_{(\bi_1^d,\ldots,\bi_{n_d}^d) \in[g(I_d)]_{g(\bv_d)}}
\!\!\!\!\mu(\bi_1^d)\cdots\mu(\bi_{n_d}^d)F(\bi_1^d,\ldots,\bi_{n_d}^d)\bigg)\nonumber
\\
&\leq
 f(\bv)^{d-1}\sum_{g\in G}\bigg(K_1 \mu(C_{g(\bv_1)})^{1/\a_1}
\prod_{i=1}^{n_1-1}
 \Big(\sum_{\bj\in L_i^1, \bj \succeq g(\bv_1)} f(\bj)^{p_i^1}\mu(C_{\bj})^{1+p_i^1}\Big)^{1/p_i^1}\bigg)\nonumber
\\
&\hspace{2cm}\times \cdots \times\bigg(K_d\, \mu(C_{g(\bv_d)})^{1/\a_d}
\prod_{i=1}^{n_d-1}
 \Big(\sum_{\bj\in L_i^d, \bj \succeq g(\bv_d)} f(\bj)^{p_i^d}\mu(C_{\bj})^{1+p_i^d}\Big)^{1/p_i^d}\bigg)\nonumber
\\
&\mbox{(on applying (\ref{intestapnp}) to each $I_j$ in the decomposition of $I$, with $K_j$ as the corresponding constant $K$)}\nonumber
\\
&=
K_1\cdots K_d\,  f(\bv)^{d-1}\sum_{g\in G}
\bigg(\Big( \mu(C_{g(\bv_1)})^{1/\a_1}\cdots\mu(C_{g(\bv_d)})^{1/\a_d}\Big)\nonumber
\\
&\hspace{1cm}\times\prod_{i=1}^{n_1-1} \Big(\sum_{\bj\in L_i^1, \bj \succeq g(\bv_1)} f(\bj)^{p_i^1}\mu(C_{\bj})^{1+p_i^1}\Big)^{1/p_i^1}
\times \cdots \times 
\prod_{i=1}^{n_d-1}
 \Big(\sum_{\bj\in L_i^d, \bj \succeq g(\bv_d)} f(\bj)^{p_i^d}\mu(C_{\bj})^{1+p_i^d}\Big)^{1/p_i^d}\bigg)\nonumber
\\
&\leq 
K_1\cdots K_d\,  f(\bv)^{d-1}
\Big(\sum_{g\in G} \mu(C_{g(\bv_1)})^{\beta/\a_1}\cdots\mu(C_{g(\bv_d)})^{\beta/\a_d}\Big)^{1/\beta}\nonumber
\\
&\hspace{1cm}\times\prod_{i=1}^{n_1-1} \Big(\sum_{g\in G}\sum_{\bj\in L_i^1, \bj \succeq g(\bv_1)} f(\bj)^{p_i^1}\mu(C_{\bj})^{1+p_i^1}\Big)^{1/p_i^1}
\times \cdots \times 
\prod_{i=1}^{n_d-1}
 \Big(\sum_{g\in G}\sum_{\bj\in L_i^d, \bj \succeq g(\bv_d)} f(\bj)^{p_i^d}\mu(C_{\bj})^{1+p_i^d}\Big)^{1/p_i^d}\nonumber
\\
&\mbox{(using H\"{o}lder's inequality, noting that $\frac{1}{\beta}+ \sum_{j=1}^{d}\sum_{i=1}^{n_j -1} \frac{1}{p_i^j} = 1$)}\nonumber
\\
&\leq
K_1\cdots K_d\,  f(\bv)^{d-1}
K(m;\beta/\a_1,\ldots,\beta/\a_d,0,\ldots, 0) ^{1/\beta} \big(\mu(C_{\bv_1})+\cdots+\mu(C_{\bv_m})\big)^{1/\a_1 + \cdots +1/\a_d}\nonumber
\\
&\hspace{1cm}\times\prod_{i=1}^{n_1-1} \Big((m-1)! \!\!\!\! \sum_{\bj\in L_i^1, \bj \succeq \bv} f(\bj)^{p_i^1}\mu(C_{\bj})^{1+p_i^1}\Big)^{1/p_i^1}
\times \cdots \times 
\prod_{i=1}^{n_d-1}
 \Big((m-1)! \!\!\!\! \sum_{\bj\in L_i^d, \bj \succeq \bv} f(\bj)^{p_i^d}\mu(C_{\bj})^{1+p_i^d}\Big)^{1/p_i^d}\nonumber
\\
&\mbox{(where $K(m; \cdots)$ is given by (\ref{muirin1})})\nonumber
\\ 
&=
K_0 \, f(\bv)^{d-1} 
\mu(C_{\bv})^{d-1 +1/p_1^0 + \cdots +1/p_{d-1}^0 + 1/\a} \!\!\!\!
\prod_{1 \leq j \leq d, 1\leq i \leq n_{j-1}} \Big( \sum_{\bj\in L_i^j, \bj \succeq \bv} f(\bj)^{p_i^j}\mu(C_{\bj})^{1+p_i^j}\Big)^{1/p_i^j}\nonumber
\\
&\mbox{(where $K_0= K_1\cdots K_d \,K(m;\beta/\a_1,\ldots,\beta/\a_d,0,\ldots, 0) ^{1/\beta}
(m-1)! ^{1-1/\beta}$ and using (\ref{betadef}))}\label{kuse}
\\
&=K_0\,
\mu(C_{\bv})^{1/\a}\big( f(\bv)^{p_1^0}\mu(C_{\bv})^{1+p_1^0}\big)^{1/p_1^0}
\cdots \big( f(\bv)^{p_{d-1}^0}\mu(C_{\bv})^{1+p_{d-1}^0}\big)^{1/p_{d-1}^0}\nonumber\\
&\hspace{1cm}\times  \prod_{1 \leq j \leq d, 1\leq i \leq n_{j-1}} \Big( \sum_{\bj\in L_i, \bj \succeq \bv} f(\bj)^{p_i^j}\mu(C_{\bj})^{1+p_i^j}\Big)^{1/p_i^j}\nonumber
\\
&=K_0\,
\mu(C_{\bv})^{1/\a}
\prod_{0 \leq j \leq d, 1\leq i \leq n_{j-1}} \Big(\sum_{\bj\in L_i, \bj \succeq \bv} f(\bj)^{p_i^j}\mu(C_{\bj})^{1+p_i^j}\Big)^{1/p_i^j}\nonumber
\\
&\mbox{(on incorporating the terms involving $p_0^j$ as single terms in the product)}\nonumber
\\
&=K_0\,
\mu(C_{\bv})^{1/\a}
\prod_{i=1}^{n-1} \Big(\sum_{\bj\in L_i, \bj \succeq \bv} f(\bj)^{p_i}\mu(C_{\bj})^{1+p_i}\Big)^{1/p_i}.\nonumber
 \end{align}
 Thus \textbf{P}(\textit{n}) holds with $K = K_0/(m-d)!$.
\medskip

\noindent{\it Case} (b2).  
Now suppose that $\bv$ is not a join point of $I$. As before let $p_1, \ldots,p_{n-1}>0$  be the exponents associated with the join points of $I$ where $\sum_{i=1}^{n-1} \frac{1}{p_i} = 1- \frac{1}{\alpha}$ with $1 <\alpha \leq \infty$.  Let $\bw$ be the first join point of $I$ below $\bv$, so that  $\bi_1,\ldots,\bi_n\succeq \bw\succeq \bv$, with $n\geq 2$.  Let $\bv_1, \ldots, \bv_r$  be the vertices of $T$ such that   $\bv_j\succeq \bv$ and $|\bv_j| = |\bw|$, and let $g_j \in  \mbox{Aut}_\bv$ be such that
$g_j(\bw) = \bv_j$.
We may decompose the orbits of $I$ as 
$$[I]_\bv = \{[g_j(I)]_{\bv_j}\}_{j=1}^r.$$
Applying case (b1) to $[g_j(I)]_{\bv_j}$ for $1 \leq j \leq r$, and letting $K$ be the corresponding constant given by (\ref{intestapnp}), which will be the same for each $j$, 
we obtain
\begin{align*}
\sum_{(\bi_1,\ldots,\bi_n) \in[I]_\bv} \mu(\bi_1)\cdots\mu(\bi_n)&F(\bi_1,\ldots,\bi_n)\\
&=
\sum_{j=1}^r \sum_{(\bi_1,\ldots,\bi_n) \in[g_j(I)]_{\bv_j}}     
\mu(\bi_1)\cdots\mu(\bi_n)F(\bi_1,\ldots,\bi_n)\\
&\leq 
\sum_{j=1}^r K \mu(C_{\bv_j})^{1/\alpha}\prod_{i=1}^{n-1}
 \Big(\sum_{\bj\in L_i, \bj \succeq \bv_j} f(\bj)^{p_i}\mu(C_{\bj})^{1+p_i}\Big)^{1/p_i}\\
 &\leq 
K\Big(\sum_{j=1}^r \mu(C_{\bv_j})\Big)^{1/\alpha}\prod_{i=1}^{n-1}
 \Big(\sum_{j=1}^r\sum_{\bj\in L_i, \bj \succeq \bv_j} f(\bj)^{p_i}\mu(C_{\bj})^{1+p_i}\Big)^{1/p_i}\\ 
 &= 
K\mu(C_{\bv})^{1/\alpha}\prod_{i=1}^{n-1}
 \Big(\sum_{\bj\in L_i, \bj \succeq \bv} f(\bj)^{p_i}\mu(C_{\bj})^{1+p_i}\Big)^{1/p_i}\ \end{align*}
using H\"{o}lder's inequality, to get (\ref{intestapnp}) in this case.

Thus by induction \textbf{P}(\textit{n}) holds for all $n \geq 1$, and the Lemma follows if $n >1$ taking $\alpha = \infty$.
$\Box$
\medskip
 
\noindent{\it Proof of Theorem \ref{thmmary}}  
The theorem is immediate on setting $\bv = \emptyset$ in (\ref{intestapnp}).
$\Box$

\section{Value of the constant $K$}\label{valk}
\setcounter{equation}{0}
\setcounter{theo}{0}

The constant $K$ that occurs in (\ref{intesta}) arises from a product of terms  $K(m;\beta/\a_1,\ldots,\beta/\a_d,0,\ldots, 0)$ that are incorporated at (\ref{kuse}) at each step of the induction.
Thus to estimate $K$ we first need to bound $K(m;a_1,\ldots,a_m)$.

\begin{lem}\label{muirhead}
Let $m \geq 1$ and let $a_i \geq 0$ for $i=1,2,\ldots,m$, with $0 <a_1 + a_2 + \cdots + a_m =  s$. Let $S_m$ be the symmetric group of all permutations $\sigma$ of $(1,2,\ldots,m)$. Recall that
$K(m;a_1, \ldots,a_m)$  is the least number such that 
\be
  \sum_{\sigma \in S_m} x_{\sigma(1)}^{a_1}x_{\sigma(2)}^{a_2}\cdots x_{\sigma(m)}^{a_m}
  \leq K(m;a_1, \ldots,a_m)(x_1 + x_2 + \cdots +x_m)^s \label{muirin}
\ee
for all $x_i \geq 0 $, with the convention that $0^0 = 1$. 

(i) If $0 <s \leq 1$ then
$$K(m;a_1, \ldots,a_m)= m! \, m^{-s}.$$

(ii) If $1\leq s$ then
$$m! \, m^{-s} \leq K(m;a_1, \ldots,a_m)\leq (m-1)!.$$

(iii)  If $1\leq s$ and $a_i \geq (s-1)/m$ for  all $i=1,\ldots,m$ then
$$K(m;a_1, \ldots,a_m)= m! \, m^{-s}.$$

(iv) For $m=2$, if $(a_1- a_2)^2 \leq a_1+ a_2 $ then
$$K(2;a_1,a_2)= 2^{1-s}=2^{1-a_1-a_2}.$$
With the values given in (i), (iii) and (iv) there is equality in (\ref{muirin})  if and only if the $x_i$ are all equal.
\end{lem}
\noindent{\it Proof.}  
It is enough to prove (\ref{muirin}) with the $x_i>0$ and take the limit of the inequality for any  $x_i = 0$.
Modifying the exponents to sum to 1 and then using Muirhead's inequality (see \cite{Gar,HLP}) gives
\begin{align*}
\frac{1}{m!}\sum_{\sigma \in S_m} x_{\sigma(1)}^{a_1}x_{\sigma(2)}^{a_2}\cdots x_{\sigma(n)}^{a_m}
&  =
\frac{1}{m!}\sum_{\sigma \in S_m} (x_{\sigma(1)}^s)^{a_1/s}(x_{\sigma(2)}^s)^{a_2/s}\cdots (x_{\sigma(m)}^s)^{a_m/s}\\
&\leq
\frac{1}{m}(x_1^s + x_2^s + \cdots +x_m^s).
\end{align*}
Then (i) follows directly on applying the generalized mean inequality and (ii) follows using Minkowski's inequality. Note that in  case  (ii), setting $K= m!m^{-s}$  in (\ref{muirin}) we get equality when the $x_i$ are all equal but the inequality may fail with this $K$ for other $x_i$. For (iii) we have
\begin{align*}
\frac{1}{m!}\sum_{\sigma \in S_m} x_{\sigma(1)}^{a_1}x_{\sigma(2)}^{a_2}\cdots x_{\sigma(n)}^{a_m}
&  =
\frac{1}{m!}\sum_{\sigma \in S_m}x_{\sigma(1)}^{a_1-(s-1)/m}\cdots x_{\sigma(n)}^{a_m-(s-1)/m} \big( x_{\sigma(1)}x_{\sigma(2)}\cdots x_{\sigma(m)}\big)^{(s-1)/m}  \\
&\leq
\frac{1}{m}(x_1 + x_2 + \cdots +x_m)\big((x_1 + x_2 + \cdots +x_m)/m\big)^{s-1}
\end{align*}
using Muirhead's inequality and the geometric-arithmetic mean inequality.

For (iv) we first claim that for all $r,q \geq 0$ such that $(r- q)^2 \leq r+q $
\be
\frac{\cosh(r-q)\theta}{(\cosh \theta)^{r+q}} < 1\quad \mbox{ for all } 0 \neq \theta \in \mathbb{R}. \label{coshrat}
\ee
To see this, assume, without loss of generality, that $0 \leq q<r $. Differentiating and using the addition formula for hyperbolic functions gives
\be
\frac{{\rm d}}{{\rm d}\theta} \bigg[\frac{\cosh(r-q)\theta}{(\cosh \theta)^{r+q}}\bigg]
=\frac{r \sinh\big((r-q-1)\theta\big) - q \sinh\big((r-q+1)\theta\big)}{(\cosh \theta)^{r+q} \cosh \theta}.\label{diff}
\ee
Note that
$$r(r-q-1) - q(r-q+1) = (r- q)^2 - (r+q) \leq 0$$
from our assumption. Since $(r-q-1)^2 < (r-q+1)^2$ and $q(r-q+1)\geq 0$ it follows inductively that
\be
r(r-q-1)^{2k+1} - q(r-q+1)^{2k+1}< 0\label{twok}
\ee
for all integers $k \geq 1$.
But the left-hand expression in (\ref{twok}) is just the coefficient of $\theta^{2k+1}/(2k+1)! $ in the power series expansion of the  numerator of  the quotient in (\ref{diff}). Thus the derivative (\ref{diff}) is strictly positive for $\theta<0$ and strictly negative for $\theta>0$, from which (\ref{coshrat}) follows.

If $x_             1= \lambda x_2$ where $\lambda >0$ and setting $\lambda^{1/2} = {\rm e}^\theta$,
$$\frac{x_1^{r}x_2^{q} +x_2^{r}x_1^{q}}{(x_1+x_2)^{r + q}}
= \frac{\lambda^{(r-q)/2}+ \lambda^{(q-r)/2} }{(\lambda^{1/2}+ \lambda^{-1/2})^{r+q} }
= \frac{{\rm e}^{(r-q)\theta}+ {\rm e}^{-(r-q)\theta} }{({\rm e}^{\theta}+{\rm e}^{-\theta})^{r+q} }
= 2^{1- r-q}\frac{\cosh(r-q)\theta}{(\cosh \theta)^{r+q}} \leq 2^{1- r-q}$$
by (\ref{coshrat}), with equality if and only if $x_1=  x_2$. Taking $r=a_1$ and $q=a_2$ gives (iv).
$\Box$
\medskip

Note that case (iv) of inequality (\ref{muirin}) when $m=2$ is related to Muirhead means and Schur convexity, which has a substantial literature, see \cite{CX,Gar}. Some related inequalities are obtained in \cite{CX, WCQ} but we were unable to find the particular inequality that we required.

Note also in relation to case (iv) that if $(a_1- a_2)^2 > a_1+ a_2$ then the maximum of 
$ \sum_{\sigma \in S_2} x_{\sigma(1)}^{a_1}x_{\sigma(2)}^{a_2}\big/
  (x_1 + x_2 )^{a_1 + a_2}$ does not necessarily occur when $x_1 = x_2$, in which case $K(2; a_1,a_2) > 2^{1-a_1-a2}$.
 \medskip

\noindent{\it Proof of Corollary \ref{genk}} 
Each time the induction step is applied at a join point with multiplicity $d-1$      the constant gets multiplied at $(\ref{kuse})$ by 
\begin{align*}
K(m;\beta/\a_1,&\ldots,\beta/\a_d,0,\ldots, 0) ^{1/\beta}(m-1)!^{1-1/\beta} (m-d)!^{-1}\\
&\leq (m-1)!^{1/\beta}(m-1)!^{1-1/\beta}  (m-d)!^{-1}= (m-1)! (m-d)!^{-1}
\end{align*}
using Lemma \ref{muirhead}(ii). 
The induction step is applied once at each join  point, so the estimate for $K$ follows. Noting that 
$(m-1)!/(m-r_i-1)! \leq (m-1)^{r_i}$ for each $i$ and $\sum_{i=1}^j r_i= n-1$ gives the stated inequality.
$\Box$
\medskip

In general the value of $K$ stated in Corollary \ref{genk} will not be optimal which in the case of a binary tree gives $K=1$. However, provided that the exponents $p_i$ are reasonably well distributed in the manner stated precisely in Corollary \ref{thmbinary}, we can reduce this to $K=2^{-(n-1)}$. 
\medskip

\noindent{\it Proof of Corollary \ref{thmbinary}} 
For a binary tree, $m=2$, all the join points are of multiplicity $1$, and the inductive step is applied $n-1$ times. Thus the corollary will follow if  we can show that the multiplier incorporated at \eqref{kuse} each time the inductive step is used satisfies 
$$K(2;\beta/\a_1,\beta/\a_2)^{1/\beta} \leq 2^{-1}.$$

Note that (\ref{halves}) implies that at each stage of the induction $1/\a_1, 1/\a_2 \geq \frac{1}{2}$ each time (\ref{ajs}) is used. In particular
$(1/\a_1 -1/4)(1/\a_2 -1/4)\geq 1/16$
 which rearranges to
\begin{equation}
\frac{1}{\a_1} + \frac{1}{\a_2} \leq \frac{4}{\a_1 \a_2}.\label{ain}
\end{equation}
Then
\begin{eqnarray*}\
\bigg(\frac{1}{\a_1} - \frac{1}{\a_2}\bigg)^2
&=&\frac{1}{\a_1^2} + \frac{1}{\a_2^2}-\frac{2}{\a_1 \a_2}\\
&\leq& \bigg(\frac{1}{\a_1^2} + \frac{1}{\a_2^2}+\frac{2}{\a_1 \a_2}\bigg)- \bigg(\frac{1}{\a_1} + \frac{1}{\a_2}\bigg)\\
&=&\bigg(\frac{1}{\a_1} + \frac{1}{\a_2}\bigg)\bigg(\frac{1}{\a_1} + \frac{1}{\a_2}-1\bigg)\\
&=&\bigg(\frac{1}{\a_1} + \frac{1}{\a_2}\bigg)\frac{1}{\beta}
\end{eqnarray*}
using  (\ref{betadef}), giving
$$\bigg(\frac{\beta}{\a_1} - \frac{\beta}{\a_2}\bigg)^2 \leq \frac{\beta}{\a_1} +\frac{\beta}{\a_2}.$$
Thus Lemma \ref{muirhead} (iv) gives
$$K(2;\beta/\a_1,\beta/\a_2)^{1/\beta} \leq 2^{(1-\beta/\a_1-\beta/\a_2)/\beta}=2^{-1}$$
as required. 
$\Box$ 
\medskip

Finally note that conditions for equality are not in general easy to specify. This requires equality at each step of the induction when 
H\"{o}lder's inequality and  (\ref{muirin}) are applied. There will be equality in Theorem \ref{thmmary}  if $\mu(\bi)$ is constant for all $\bi \in T_0$ and $f(\bj)$ is constant for all $\bj \in L_i$ for each level $L_i$ provided that we have the optimal value of $K$. This we can do under the conditions of Corollary \ref{thmbinary}. However, as the value of $K(m;a_1, \ldots,a_m)$ for Lemma \ref{muirhead} case (iii) seems to be unknown when $m\geq 3$ and $s>1$,  the value of $K$ given by Proposition \ref{genk} will not in general be optimal. 

\section{Acknowledgement}
\setcounter{equation}{0}
\setcounter{theo}{0}
The author is most grateful to the referee for a number of helpful suggestions.

\end{document}